\documentclass[3p,times]{elsarticle}
\usepackage{amsfonts}
\usepackage{mathrsfs}
\usepackage{amsmath}
\usepackage{amssymb}

 \usepackage{lineno}
\newtheorem{theorem}{\indent Theorem}[section]

\journal{Journal of Differential Equations}

\begin{document}

\begin{frontmatter}



\title{Global attractor of the three dimensional primitive equations of large-scale ocean and atmosphere dynamics in an unbounded domain}


\author{Bo You$^{a,*}$, Fang Li$^b$}

\address{$^a$ School of Mathematics and Statistics, Xi'an Jiaotong University, Xi'an, 710049, P. R. China\\
$^b$ School of Mathematics and Statistics, Xidian University, Xi'an, 710126, P. R. China}

 \cortext[cor1]{Corresponding author.\\
 E-mail addresses:youb2013@xjtu.edu.cn(B. You),fli@xidian.edu.cn (F. Li)}
%\author{Bo You\corref{cor1}}
%\ead{youb2013@mail.xjtu.edu.cn}
% \address{School of Mathematics and Statistics, Xi'an Jiaotong University,
% Xi'an, 710049, P. R. China}
% \author{Fang Li\corref{cor2}}
%\ead{lifang101216@126.com}
% \address{Department of Mathematics, Nanjing University\\
%    Nanjing, 210093, P. R. China}
%\author{Chengkui Zhong\corref{cor2}}
%\ead{ckzhong@nju.edu.cn}
% \address{Department of Mathematics, Nanjing University\\
%    Nanjing, 210093, P. R. China}

\begin{abstract}
This paper is concerned with the long-time behavior of solutions for the three dimensional primitive equations of large-scale ocean and atmosphere dynamics in an unbounded domain. Since the Sobolev embedding is no longer compact in an unbounded domain, we cannot obtain the asymptotical compactness of the semigroup generated by problem \eqref{2.4}-\eqref{2.6} by the Sobolev compactness embedding theorem even if we obtain
the existence of an absorbing set in more regular phase space. To overcome this difficulty, we first prove the asymptotical compactness in a weak phase space of the semigroup generated by problem \eqref{2.4}-\eqref{2.6} by combining the tail-estimates and the energy methods. Thanks to the existence of an absorbing set in more regular phase space, we establish the existence of a global attractor of problem \eqref{2.4}-\eqref{2.6} by virtue of Sobolev interpolation inequality.
\end{abstract}


\end{frontmatter}

%%\linenumbers
%%\begin{linenumbers}

\section{Introduction}
\def\theequation{1.\arabic{equation}}\makeatother
\setcounter{equation}{0}
In this paper, we consider the long-time behavior of solutions for the following three dimensional primitive equations of large-scale ocean and atmosphere dynamics(see \cite{pj1})
\begin{equation}\label{1.1}
\begin{cases}
&\frac{\partial v}{\partial t}+(v\cdot \nabla)v+w\frac{\partial
v}{\partial z}+\nabla
p+\frac{1}{Ro}f{\vec{k}}\times v+L_1v=0,\\
&\frac{\partial p}{\partial z}+T=0,\\
&\nabla\cdot v+\frac{\partial w}{\partial z}=0,\\
&\frac{\partial T}{\partial t}+v\cdot \nabla
T+w\frac{\partial T}{\partial z}+L_2T=Q
\end{cases}
\end{equation}
in an unbounded domain
\begin{align*}
\Omega=\mathbb{R}\times(0,l)\times(-h,0),
\end{align*}
where $l,$ $h$ are two positive constants
and the horizontal velocity field $v=(v_1,v_2)$, the
three-dimensional velocity field ($v_1,$ $v_2,$ $w$), the
temperature $T$ and the pressure $p$ are the unknowns,  $\vec{k}=(0,0,1)$ is
vertical unit vector in $\mathbb{R}^3$, $f=f_0+\beta y$ is the Coriolis parameter,
$Ro$ is the Rossby number, which measures the significant
influence of the rotation of the earth to the dynamical behaviour of the
ocean, $Q(x,y,z)$ is a given heat source. The viscosity and the
heat diffusion operators $L_1$ and $L_2$ are given by
\begin{align*}
&L_1=-\frac{1}{Re_1}\Delta-\frac{1}{Re_2}\frac{\partial^2}{\partial
z^2},\\
&L_2=-\frac{1}{Rt_1}\Delta-\frac{1}{Rt_2}\frac{\partial^2}{\partial
z^2},
\end{align*}
where $Re_1,$ $Re_2$ are positive constants representing
the horizontal and vertical Reynolds numbers, respectively, and $
Rt_1,$ $
Rt_2$ are positive constants which stand for the
horizontal and vertical eddy diffusivity, respectively. For the sake of simplicity, let $\nabla=(\partial_x,\partial_y)$ be the horizontal gradient
operator and let $\Delta=\partial_x^2+\partial_y^2$ be
the horizontal Laplacian. We observe that the above system is
similar to the three dimensional Boussinesq system where the equation of
vertical motion is approximated by the hydrostatic balance.

Let $\Gamma_u$ and $\Gamma_b$ be the upper and bottom boundaries of $\Omega$ given by
\begin{align*}
&\Gamma_u=\{(x,y,z)\in \overline{\Omega}:z=0\},\\
&\Gamma_b=\{(x,y,z)\in \overline{\Omega}:z=-h\}.
\end{align*}
 Then we can state the boundary conditions of the system \eqref{1.1}
as follows
\begin{equation}\label{1.2}
\begin{cases}
& \frac{\partial v}{\partial
z}|_{\Gamma_u}=0,w|_{\Gamma_u}=0,(\frac{1}{Rt_{2}}\frac{\partial T}{\partial z}+\alpha
T)|_{\Gamma_u}=0,\\
&\frac{\partial
v}{\partial z}|_{\Gamma_b}=0,w|_{\Gamma_b}=0,\frac{\partial T}{\partial z}|_{\Gamma_b}=0,\\
&v(x,0,z,t)=v(x,l,z,t)=0,\,\frac{\partial T}{\partial y}(x,0,z,t)=\frac{\partial T}{\partial y}(x,l,z,t)=0,
\end{cases}
\end{equation}
where $\alpha$ is a positive constant related with the
turbulent heating on the surface of the ocean.

 In addition, we supply the system \eqref{1.1}-\eqref{1.2} with
the initial conditions
\begin{equation}\label{1.3}
\begin{cases}
&v(x,y,z,0)=v_{0}(x,y,z),\\
&T(x,y,z,0)=T_{0}(x,y,z).
\end{cases}
\end{equation}

Large-scale dynamics of ocean and atmosphere is governed by the
 primitive equations which are derived from the Navier-Stokes
 equations with rotation coupled to thermodynamics and salinity
 diffusion-transport equations, which account for the buoyancy
 forces and stratification effects under the Boussinesq
 approximation. Moreover, due to the shallowness of the oceans
 and the atmosphere, i.e., the depth of the fluid layer is very
 small in comparison to the radius of the earth, the vertical large-scale motion in the oceans and the atmosphere is much smaller than
 the horizontal one, which in turn leads to modeling the vertical
 motion by the hydrostatic balance. As a result, one can obtain the
 system \eqref{1.1}-\eqref{1.3} which is known as the primitive
 equations for ocean and atmosphere dynamics (see \cite{ccs1, ljl, ljl1,
pj1, vgk}). We observe that one has
 to add the diffusion-transport equation of the salinity to the
 system \eqref{1.1}-\eqref{1.3} in the case of ocean dynamics, but we omitted it here in order to simplify
 our mathematical presentation. However, we emphasize that our
 results are equally valid when the salinity effects are taken into
 account.

 In the past several decades, the well-posedness and the long-time behavior of solutions for the primitive
equations of the atmosphere, the ocean and the coupled
atmosphere-ocean have been extensively studied both from the theoretical point of view (see \cite{ccs1, da, elc, ebd, gfg, gbl2, gbl, gbl4, gbl1, gbl3, hcb1, hcb2, hcb, jn, jn2, jn1, ki, ljl2, ljl, ljl1, mtt, mtt4, mtt2, pm}) and numerical point of view (see \cite{ccs5, gfg1, hyn, hyn2, hyn1, mtt1, mtt3, sr, sj}). In particular,
 the existence and uniqueness of global
strong solutions for the
three-dimensional viscous primitive equations of large-scale ocean
was established in \cite{ccs1}. In \cite{gbl2}, the authors proved the existence of
weak solutions and trajectory attractors for the moist atmospheric
equations in geophysics. The long-time dynamics of the primitive equations of large-scale atmosphere was considered and a weakly compact global attractor
$\mathcal{A}$ which captures all the trajectories was obtained in \cite{gbl}. Under the assumption of
the initial data $(v_0,T_0)\in V\cap (H^2(\Omega))^3,$
the authors have proved the existence of a compact global attractor in $V$ for the primitive equations of large-scale atmosphere in
\cite{gbl1}. In \cite{mtt4, mtt2}, the author has proved the existence of an uniform attractor for the three dimensional non-autonomous primitive equations with oscillating external force and investigated its structure and the upper semi-continuous properties. The global well-posedness and long-time dynamics
for the three-dimensional stochastic primitive equations of large-scale ocean was considered and the existence of random attractors for the corresponding random dynamical system was proved in \cite{gbl4}. In \cite{ljl2}, the authors have proved the upper bound of the fractal dimension of the global attractor for the primitive equations of atmospheric circulation and provided its physical interpretation. The existence of the global attractor for
the three dimensional viscous primitive equations of large-scale atmosphere in log-pressure coordinate was
established in \cite{yb2}. The existence of a global attractor in $V$ for the primitive equations of large-scale atmosphere and ocean dynamics was proved by Ning Ju in \cite{jn} by using the Aubin-Lions compactness theorem under the assumption $Q\in L^2(\Omega).$ In \cite{jn2, jn1}, the authors have proved the finite dimensional global attractor for the 3D viscous primitive equations by using the squeezing property. However, there are no results related to the existence of a global attractor for the three dimensional primitive equations of large-scale ocean and atmosphere dynamics in an unbounded domain.

The main purpose of this paper is to study the existence of a global attractor in $V$ for the three dimensional primitive equations of large-scale ocean and atmosphere dynamics in an unbounded domain with $Q\in H^1(\Omega).$ First of all, we make some a priori estimates for the solutions of problem \eqref{2.4}-\eqref{2.6}, which imply the existence of an absorbing set in $(H^2(\Omega))^3\cap V$ for the semigroup $\{S(t)\}_{t\geq0}$ generated by problem \eqref{2.4}-\eqref{2.6}. Since the Sobolev embedding is no longer compact in an unbounded domain, we cannot prove the asymptotical compactness in $V$ of the semigroup $\{S(t)\}_{t\geq0}$ generated by problem \eqref{2.4}-\eqref{2.6} by the Sobolev compactness embedding theorem even if we can establish
the existence of an absorbing set in more regular phase space than $V.$ To overcome this difficulty, we first prove the asymptotical compactness in $H$ of the semigroup $\{S(t)\}_{t\geq0}$ generated by problem \eqref{2.4}-\eqref{2.6} by combining the tail-estimates for the equation of the temperature and the energy methods for the equation of the velocity. Finally, we prove the existence of a global attractor in $V$ for the semigroup $\{S(t)\}_{t\geq0}$ generated by problem \eqref{2.4}-\eqref{2.6} via Sobolev interpolation inequality.

 Throughout this paper, let $X$ be a Banach space endowed with the norm $\|\cdot\|_X,$ denote by $\|u\|_p$ the
  $L^p(\Omega)$-norm of $u$ and let $C$ be the positive generic constants independent of the initial data,
  which may be different from line to line.

\section{New formulation and functional setting}
\def\theequation{2.\arabic{equation}}\makeatother
\setcounter{equation}{0}
\subsection{New formulation}
Integrating the third equation of \eqref{1.1} in the $z$ direction, we obtain
\begin{align*}
w(x,y,z,t)=w(x,y,-h,t)-\int_{-h}^{z}\nabla\cdot v(x,y,r,t)\,dr.
\end{align*}

    Since $w(x,y,0,t)=w(x,y,-h,t)=0$ by \eqref{1.2}, we have
\begin{align}\label{2.1}
w(x,y,z,t)=-\int_{-h}^{z}\nabla\cdot v(x,y,r,t)\,dr
\end{align}
and
\begin{align}\label{2.2}
\int_{-h}^{0}\nabla\cdot
v(x,y,r,t)\,dr=\nabla\cdot\int_{-h}^{0}
v(x,y,r,t)\,dr=0.
\end{align}

Integrating the second equation of \eqref{1.1} with respect to $z,$ we obtain
\begin{align}\label{2.3}
p(x,y,z,t)=p_s(x,y,t)-\int_0^{z}T(x,y,r,t)\,dr,
\end{align}
where $p_s(x,y,t)$ is a free function to be determined.

    Based on \eqref{2.1} and \eqref{2.3}, we obtain the following new formulation of \eqref{1.1}-\eqref{1.3}:
\begin{equation}\label{2.4}
\begin{cases}
&\frac{\partial v}{\partial t}+L_1v+(v\cdot
\nabla)v-\left(\int_{-h}^z\nabla\cdot
v(x,y,r,t)\,dr\right)\frac{\partial v}{\partial z}+\nabla
p_s(x,y,t)+\frac{1}{Ro}f{\vec{k}}\times v=\int_0^{z}\nabla T(x,y,r,t)\,dr,\\
& \frac{\partial T}{\partial t}+L_2T+v\cdot\nabla
T-\left(\int_{-h}^z\nabla\cdot v(x,y,r,t)\,dr\right)\frac{\partial
T}{\partial z}=Q
\end{cases}
\end{equation}
with the following boundary conditions
\begin{equation}\label{2.5}
\begin{cases}
&\frac{\partial v}{\partial z}|_{\Gamma_u}=0,\frac{\partial
v}{\partial z}|_{\Gamma_b}=0,v(x,0,z,t)=v(x,l,z,t)=0,\\
&(\frac{1}{Rt_2}\frac{\partial T}{\partial z}+\alpha
T)|_{\Gamma_u}=0,\frac{\partial T}{\partial z}|_{\Gamma_b}=0,
\,\frac{\partial T}{\partial y}(x,0,z,t)=\frac{\partial T}{\partial y}(x,l,z,t)=0
\end{cases}
\end{equation}
and the initial data
\begin{equation}\label{2.6}
\begin{cases}
&v(x,y,z,0)=v_0(x,y,z),\\
&T(x,y,z,0)=T_0(x,y,z).
\end{cases}
\end{equation}

 Denote
\begin{align*}
 \bar{v}(x,y)=\frac{1}{h}\int_{-h}^0v(x,y,r)\,dr
\end{align*}
and
\begin{align*}
  \tilde{v}=v-\bar{v}.
\end{align*}

Taking the average of the first equation of \eqref{2.4} in the $z$-direction over the interval $(-h,0)$ and using \eqref{2.5}, we obtain
\begin{align}\label{2.7}
\frac{\partial \bar{v}}{\partial
t}-\frac{1}{Re_1}\Delta\bar{v}+(\bar{v}\cdot\nabla)\bar{v}+\overline{(\tilde{v}\cdot\nabla)\tilde{v}+(\nabla\cdot
\tilde{v})\tilde{v}}+\nabla
p_{s}(x,y,t)+\frac{1}{Ro}f{\vec{k}}\times
\bar{v}-\overline{\int_0^z\nabla T(x,y,r,t)\,dr} =0
\end{align}
with the boundary conditions
\begin{align}\label{2.8}
\nabla\cdot \bar{v}=0, \bar{v}(x,0,t)=\bar{v}(x,l,t)=0.
\end{align}
 Subtracting \eqref{2.7} from the first equation of \eqref{2.4}, we obtain
\begin{align}\label{2.9}
\nonumber&\frac{\partial \tilde{v}}{\partial t}+L_1\tilde{v}+(\tilde{v}\cdot
\nabla)\tilde{v}-\left(\int_{-h}^{z}\nabla\cdot
\tilde{v}(x,y,r,t)\,dr\right)\frac{\partial \tilde{v}}{\partial
z}+(\tilde{v}\cdot\nabla)\bar{v}+(\bar{v}\cdot\nabla)\tilde{v}-\int_0^z\nabla T(x,y,r,t)\,dr\\
&+\frac{1}{Ro}f{\vec{k}}\times \tilde{v}-\overline{(\tilde{v}\cdot\nabla)\tilde{v}+(\nabla\cdot
\tilde{v})\tilde{v}}+\overline{\int_0^{z}\nabla T(x,y,r,t)\,dr
}=0
\end{align}
with the boundary conditions
 \begin{align}\label{2.10}
\frac{\partial \tilde{v}}{\partial
z}|_{\Gamma_u}=0,\frac{\partial \tilde{v}}{\partial
z}|_{\Gamma_b}=0,\tilde{v}(x,0,z,t)=\tilde{v}(x,l,z,t)=0.
\end{align}
\subsection{Functional spaces and some lemmas}
To study problem \eqref{2.4}-\eqref{2.6}, we introduce some functional spaces. Define
\begin{align*}
\mathcal {V}_1=&\left\{v\in (C^{\infty}(\bar{\Omega}))^{2}:\frac{\partial
v}{\partial z}|_{\Gamma_u}=0,\frac{\partial v}{\partial
z}|_{\Gamma_b}=0,v(x,0,z)=v(x,l,z)=0,\int_{-h}^{0}\nabla\cdot v(x,y,r)\,dr=0\right\},\\
\mathcal {V}_2=&\left\{T\in
C^{\infty}(\bar{\Omega}):(\frac{1}{Rt_{2}}\frac{\partial T}{\partial
z}+\alpha T)|_{\Gamma_u}=0,\frac{\partial T}{\partial z}|_{\Gamma_b}=0,
\frac{\partial T}{\partial y}(x,0,z)=\frac{\partial T}{\partial y}(x,l,z)=0\right\}.
\end{align*}
Denote by $ V_1$ and $V_2,$ respectively, be the closure of $\mathcal {V}_1$ and $\mathcal {V}_2$ with respect to the following
norms given by:
\begin{align*}
\|v\|^2=&\frac{1}{Re_1}\int_{\Omega}|\nabla v|^2\,dxdydz+\frac{1}{Re_2}\int_{\Omega}|\partial_z v|^2\,dxdydz,\\
\|T\|^2=&\frac{1}{Rt_1}\int_{\Omega}|\nabla T|^2\,dxdydz+\frac{1}{Rt_2}\int_{\Omega}|T_z|^2\,dxdydz+\alpha\int_{\mathbb{R}}\int_0^l|T(x,y,0)|^2\,dxdy
\end{align*}
for any $v\in\mathcal {V}_1$ and any $T\in\mathcal {V}_2,$  let $ H_1$ be the closure of $ \mathcal{V}_1$ with respect to the $(L^2(\Omega))^2$-norm, let
$V=V_1\times V_2$ and $H=H_1\times L^2(\Omega),$ denote by $X'$ the dual space of the Banach space $X.$

\section{Some a priori estimates of strong solutions}
\def\theequation{3.\arabic{equation}}\makeatother
\setcounter{equation}{0}
\subsection{The well-posedness of strong solutions}
The well-posedness of strong solutions for the three dimensional viscous primitive equations of large-scale ocean and atmosphere dynamics \eqref{2.4}-\eqref{2.6} can be obtained by the Faedo-Galerkin methods (see \cite{tr}). Now, we only state it as follows.
\begin{theorem}\label{3.1.1}
Assume that $ Q\in L^2(\Omega) $. Then for each $ (v_0, T_0)\in V,$ there exists a unique
strong solution $(v,T)\in \mathcal{C}(\mathbb{R}^+;V)$ of problem \eqref{2.4}-\eqref{2.6}, which depends continuously on the initial data in $V.$
\end{theorem}

By Theorem \ref{3.1.1}, we can define the operator semigroup $\{S(t)\}_{t\geq 0}$ in $V$ by
\begin{align*}
S(t)(v_0,T_0)=(v(t),T(t))=(v(t;(v_0,T_0)),T(t;(v_0,T_0)))
\end{align*}
for all $t\geq0,$ which is $(V,V)$-continuous, where $(v(t),T(t))$ is the strong solution of problem \eqref{2.4}-\eqref{2.6} with $(v(x,0),T(x,0))=(v_0,T_0)\in V.$
\subsection{Some a priori estimates of strong solutions}
In this subsection, we give some a priori estimates of strong solutions for problem \eqref{2.4}-\eqref{2.6}, which imply the existence of an absorbing set in $(H^2(\Omega))^3\cap V$ for the semigroup $\{S(t)\}_{t\geq 0}$ generated by problem \eqref{2.4}-\eqref{2.6}.
\subsubsection{$L^2(\Omega)$ estimates of $T$}
 Taking the inner product of the second equation of \eqref{2.4} with $T$ in $L^2(\Omega),$ we obtain
\begin{align}\label{3.2.1}
\frac{1}{2}\frac{d}{dt}\|T(t)\|^2_2+\|T(t)\|^2=\int_{\Omega}QT\,dxdydz.
\end{align}
Thanks to
\begin{align*}
\|T(x,y,z)\|^2_2\leq 2h\int_{\mathbb{R}}\int_0^l|T(x,y,0)|^2\,dxdy+2h^2\|
T_z(x,y,z)\|^2_2,
\end{align*}
we find
\begin{align}\label{3.2.2}
\frac{\|T\|^2_2}{2Rt_2h^2+\frac{2h}{\alpha}}
\leq  \frac{1}{Rt_2}\int_{\Omega}|
T_z|^2\,dxdydz+\alpha\int_{\mathbb{R}}\int_0^l|T(x,y,0)|^2\,dxdy\leq\|T\|^2.
\end{align}
It follows from \eqref{3.2.1}-\eqref{3.2.2} that
\begin{align*}
\frac{d}{dt}\|T(t)\|^2_2+\|T(t)\|^2\leq(2Rt_2h^{2}+\frac{2h}{\alpha})\|Q\|^2_2.
\end{align*}
By using \eqref{3.2.2} again, we obtain
\begin{align*}
\frac{d}{dt}\|T(t)\|^2_2+\frac{\|
T(t)\|^2_2}{2Rt_2h^2+\frac{2h}{\alpha}} \leq
(2Rt_2h^2+\frac{2h}{\alpha})\|Q\|^2_2.
\end{align*}
From the classical Gronwall inequality, we deduce
\begin{align*}
\|T(t)\|_2^2\leq \|T_0\|_2^2\exp(\frac{-t}{2Rt_2h^2+\frac{2h}{\alpha}})+(2Rt_2h^2+\frac{2h}{\alpha})^2\|Q\|_2^2,
\end{align*}
which implies that there exists a positive constant $\rho_1$ and some time $T_1>0$ such that
\begin{align}\label{3.2.3}
\|T(t)\|_2^2+\int_{t}^{t+1}\|T(r)\|^2dr\leq\rho_1
\end{align}
for any $t\geq T_1$. For brevity, we omit writing out explicitly these bounds here and we also omit writing out other
similar bounds in our future discussion for all other uniform a priori estimates. In what follows, let $\rho_i$ for any $i\in\mathbb{Z}^+$  be the positive generic constants independent of the initial data.
\subsubsection{$H_1$ estimates of $v$}
 Multiplying the first equation of \eqref{2.4} by $v$ and integrating over $\Omega,$ we obtain
\begin{align*}
\frac{1}{2}\frac{d}{dt}\|v(t)\|^2_2+\|v(t)\|^2
\leq& h\|T\|_2\|\nabla v\|_2\\
\leq&\frac{1}{2}\|v\|^2+\frac{Re_1h^2\| T\|^2_2}{2},
\end{align*}
which entails that
\begin{align*}
\frac{d}{dt}\|v(t)\|^2_2+\|v(t)\|^2\leq Re_1h^2\|T\|^2_2.
\end{align*}
It follows from the Poinc\'{a}re inequality $\|v\|_2\leq 2l\|
\nabla v\|_2$ that
\begin{align*}
\frac{d}{dt}\|v(t)\|_2^2+\frac{1}{2lRe_1}\|v(t)\|^2_2\leq Re_1h^2\|T\|^2_2.
\end{align*}
Combining the Classical Gronwall inequality with \eqref{3.2.3}, we find
\begin{align}\label{3.2.4}
\|v(t)\|_2^2+\int_{t}^{t+1}\|v(r)\|^2\,dr\leq\rho_2
\end{align}
for any $t\geq T_2\geq T_1$.
\subsubsection{$L^6(\Omega)$ estimates of $T$}
Taking the $L^2(\Omega)$ inner product of the second equation of \eqref{2.4} with $|T|^4T,$ we have
\begin{align*}
\frac{1}{6}\frac{d}{dt}\|T(t)\|^6_6+\frac{5}{9}\||T(t)|^3\|^2\leq& \||T|^3\|^{\frac{5}{3}}_{\frac{10}{3}}\|Q\|_2\\
\leq& C\||T|^3\|^{\frac{2}{3}}_2\||T|^3\|\|Q\|_2,
\end{align*}
which implies that
\begin{align}\label{3.2.5}
\frac{d}{dt}\|T(t)\|^2_6\leq C\|Q\|_2^2.
\end{align}
Therefore, from the uniform Gronwall inequality, \eqref{3.2.5} and the Sobolev embedding Theorem, we deduce
\begin{align}\label{3.2.6}
\|T(t)\|_6\leq\rho_3
\end{align}
for any $t\geq T_2+1.$
\subsubsection{$(L^6(\Omega))^{2}$ estimates of $\tilde{v}$}
Multiplying \eqref{2.9} by $|\tilde{v}|^4\tilde{v}$ and integrating over $\Omega,$ we deduce
\begin{align}\label{3.2.7}
\nonumber&\frac{1}{6}\frac{d}{dt}\|
\tilde{v}(t)\|^6_6+\frac{1}{Re_1}\||\nabla
\tilde{v}||\tilde{v}|^2\|_2^2+\frac{1}{Re_2}\||\partial_z
\tilde{v}||\tilde{v}|^2\|_2^2+\frac{4}{9}\|
|\tilde{v}|^3\|^2\\
\nonumber\leq& C\int_{\Omega}|\bar{v}||\nabla
\tilde{v}||\tilde{v}|^5\,dxdydz+C\int_{\mathbb{R}}\int_0^l(\int_{-h}^0|T|\,dz)(\int_{-h}^0|\nabla
\tilde{v}||\tilde{v}|^4\,dz)\,dxdy\\
&+C\int_{\mathbb{R}}\int_0^l(\int_{-h}^0|\tilde{v}|^2\,dz)(\int_{-h}^0|\nabla
\tilde{v}||\tilde{v}|^4\,dz)\,dxdy.
\end{align}
It follows from H\"{o}lder inequality and Sobolev embedding theorem that
\begin{align}\label{3.2.8}
\nonumber\int_{\Omega}|\bar{v}||\nabla
\tilde{v}||\tilde{v}|^5\,dxdydz
\leq&\int_{\mathbb{R}}\int_0^l|\bar{v}|(\int_{-h}^0|\nabla
\tilde{v}|^2|\tilde{v}|^4\,dz)^{\frac{1}{2}}(\int_{-h}^0|
\tilde{v}|^6\,dz)^{\frac{1}{2}}\,dxdy\\
\leq &(\int_{\mathbb{R}}\int_0^l|\bar{v}|^4\,dxdy)^{\frac{1}{4}}\||\nabla
\tilde{v}||\tilde{v}|^2\|_2(\int_{-h}^0(\int_{\mathbb{R}}\int_0^l|\tilde{v}|^{12}\,dxdy)^{\frac{1}{2}}\,dz)^{\frac{1}{2}}
\end{align}
and
\begin{align*}
\int_{\mathbb{R}}\int_0^l|\tilde{v}|^{12}\,dxdy=\int_{\mathbb{R}}\int_0^l||\tilde{v}|^3|^4\,dxdy\leq
C\int_{\mathbb{R}}\int_0^l|\tilde{v}|^6\,dxdy\int_{\mathbb{R}}\int_0^l|\nabla|\tilde{v}|^3|^2\,dxdy,
\end{align*}
which implies that
\begin{align}\label{3.2.9}
(\int_{-h}^0(\int_{\mathbb{R}}\int_0^l|\tilde{v}|^{12}\,dxdy)^{\frac{1}{2}}\,dz)^{\frac{1}{2}}
\leq
C(\int_{\Omega}|\tilde{v}|^6\,dxdydz)^{\frac{1}{4}}(\int_{\Omega}|\nabla|\tilde{v}|^3|^2\,dxdydz)^{\frac{1}{4}}.
\end{align}
We deduce from \eqref{3.2.8}-\eqref{3.2.9} that
\begin{align}\label{3.2.10}
\int_{\Omega}|\bar{v}||\nabla \tilde{v}||\tilde{v}|^5\,dxdydz
\leq
C\|\tilde{v}\|_6^{\frac{3}{2}}\|v\|_2^{\frac{1}{2}}\|\nabla v\|_2^{\frac{1}{2}}(\int_{\Omega}|\nabla|\tilde{v}|^3|^2\,dxdydz)^{\frac{1}{4}}
(\int_{\Omega}|\nabla
\tilde{v}|^{2}|\tilde{v}|^4\,dxdydz)^{\frac{1}{2}}.
\end{align}
Similarly, we find
\begin{align}\label{3.2.11}
\nonumber\int_{\mathbb{R}}\int_0^l(\int_{-h}^0|T|\,dz)(\int_{-h}^0|\nabla
\tilde{v}||\tilde{v}|^4\,dz)\,dxdy
\leq &\int_{\mathbb{R}}\int_0^l(\int_{-h}^0|T|\,dz)(\int_{-h}^0|\nabla
\tilde{v}|^2|\tilde{v}|^4\,dz)^{\frac{1}{2}}(\int_{-h}^0|\tilde{v}|^4\,dz)^{\frac{1}{2}}\,dxdy\\
\nonumber\leq &\|T\|_6(\int_{\Omega}|\nabla
\tilde{v}|^2|\tilde{v}|^4\,dxdydz)^{\frac{1}{2}}\left(\int_{\mathbb{R}}\int_0^l(\int_{-h}^0|\tilde{v}|^4\,dz)^{\frac{3}{2}}\,dxdy\right)^{\frac{1}{3}}\\
\nonumber\leq &\|T\|_6(\int_{\Omega}|\nabla
\tilde{v}|^2|\tilde{v}|^4\,dxdydz)^{\frac{1}{2}}\left(\int_{-h}^0(\int_{\mathbb{R}}\int_0^l|\tilde{v}|^6\,dxdy)^{\frac{2}{3}}\,dz\right)^{\frac{1}{2}}\\
\leq& C\|T\|_6(\int_{\Omega}|\nabla
\tilde{v}|^2|\tilde{v}|^4\,dxdydz)^{\frac{1}{2}}\|\tilde{v}\|_6^2
\end{align}
and
\begin{align}\label{3.2.12}
\nonumber\int_{\mathbb{R}}\int_0^l(\int_{-h}^0|\tilde{v}|^{2}\,dz)(\int_{-h}^0|\nabla
\tilde{v}||\tilde{v}|^4\,dz)\,dxdy\leq &C\int_{\mathbb{R}}\int_0^l(\int_{-h}^0|\nabla
\tilde{v}|^2|\tilde{v}|^4\,dz)^{\frac{1}{2}}\int_{-h}^0|\tilde{v}|^4\,dz\,dxdy\\
\nonumber\leq &C\left(\int_{\mathbb{R}}\int_0^l(\int_{-h}^0|\tilde{v}|^4\,dz)^2\,dxdy\right)^{\frac{1}{2}}(\int_{\Omega}|\nabla
\tilde{v}|^2|\tilde{v}|^4\,dxdydz)^{\frac{1}{2}}\\
\nonumber\leq &C\int_{-h}^0(\int_{\mathbb{R}}\int_0^l|\tilde{v}|^8\,dxdy)^{\frac{1}{2}}\,dz(\int_{\Omega}|\nabla
\tilde{v}|^2|\tilde{v}|^4\,dxdydz)^{\frac{1}{2}}\\
\leq &C(\int_{\Omega}|\nabla
\tilde{v}|^2|\tilde{v}|^4\,dxdydz)^{\frac{1}{2}}\|
\tilde{v}\|^3_6\|\tilde{v}\|.
\end{align}
It follows from \eqref{3.2.7}, \eqref{3.2.10}-\eqref{3.2.12} that
\begin{align*}
\frac{d}{dt}\|
\tilde{v}(t)\|^6_6+\frac{2}{Re_1}\||\nabla
\tilde{v}||\tilde{v}|^2\|_2^2+\frac{2}{Re_2}\||\partial_z
\tilde{v}||\tilde{v}|^2\|_2^2
+2\||\tilde{v}|^3\|^2
\leq C(\|v\|^2_2\|
\nabla v\|^2_2+\|
\tilde{v}\|^2)\|
\tilde{v}\|^6_6+C\|T\|^2_6\|
\tilde{v}\|^4_6,
\end{align*}
which implies that
\begin{align}\label{3.2.13}
\frac{d}{dt}\|\tilde{v}(t)\|^2_6\leq C(\|v\|^2_2\|\nabla v\|^2_2+\|\tilde{v}\|^2)\|\tilde{v}\|^2_6+C\|T\|^2_6.
\end{align}
Therefore, we infer from \eqref{3.2.3}-\eqref{3.2.4}, \eqref{3.2.13} and the uniform Gronwall inequality that
\begin{align}\label{3.2.14}
\|\tilde{v}(t)\|^2_6+\int_t^{t+1}\||\nabla
\tilde{v}(x,y,z,r)||\tilde{v}(x,y,z,r)|^2\|_2^2\,dr\leq\rho_4
\end{align}
for any $t\geq T_2+2.$
\subsubsection{ estimates of $\int_{\mathbb{R}}\int_0^l|
\nabla\bar{v}(x,y)|^2\,dxdy$}
 Taking the $L^{2}(\mathbb{R}\times (0,l))$ inner product of equation \eqref{2.7} with
$-\Delta \bar{v}$ and combining the boundary conditions
\eqref{2.8}, we conclude
\begin{align}\label{3.2.15}
\nonumber&\frac{1}{2}\frac{d}{dt}\int_{\mathbb{R}}\int_0^l|
\nabla\bar{v}(x,y,t)|^2\,dxdy+\frac{1}{Re_1}\int_{\mathbb{R}}\int_0^l|\Delta
\bar{v}(x,y,t)|^2\,dxdy\\
\nonumber\leq &C\int_{\mathbb{R}}\int_0^l|\bar{v}||\nabla \bar{v}||\Delta
\bar{v}|\,dxdy +C\int_{\mathbb{R}}\int_0^l(\int_{-h}^0|\nabla
\tilde{v}||\tilde{v}|\,dz)|\Delta \bar{v}|\,dxdy+C\int_{\mathbb{R}}\int_0^l|\bar{v}||\Delta
\bar{v}|\,dxdy\\
&+C\int_{\mathbb{R}}\int_0^l(\int_{-h}^0|\nabla T|\,dz)|\Delta \bar{v}|\,dxdy.
\end{align}
In the following, we will estimate each term of the right hand side of \eqref{3.2.15}.
\begin{align}
\nonumber\label{3.2.16}\int_{\mathbb{R}}\int_0^l|\bar{v}||\nabla \bar{v}||\Delta\bar{v}|\,dxdy\leq&C(\int_{\mathbb{R}}\int_0^l|\bar{v}|^4\,dxdy)^{\frac{1}{4}}(\int_{\mathbb{R}}\int_0^l|\nabla
\bar{v}|^4\,dxdy)^{\frac{1}{4}}(\int_{\mathbb{R}}\int_0^l|\Delta
\bar{v}|^2\,dxdy)^{\frac{1}{2}}\\
\leq &C(\int_{\mathbb{R}}\int_0^l|
\bar{v}|^2\,dxdy)^{\frac{1}{4}}
(\int_{\mathbb{R}}\int_0^l|\nabla\bar{v}|^2\,dxdy)^{\frac{1}{2}}(\int_{\mathbb{R}}\int_0^l|\Delta
\bar{v}|^2\,dxdy)^{\frac{3}{4}},\\
\label{3.2.17}\int_{\mathbb{R}}\int_0^l(\int_{-h}^{0}|\tilde{v}||\nabla
\tilde{v}|\,dz)|\Delta \bar{v}|\,dxdy
\leq &C(\int_{\Omega}|\nabla
\tilde{v}|^2|\tilde{v}|^4\,dxdydz)^{\frac{1}{4}}\|\nabla\tilde{v}\|_2^{\frac{1}{2}}(\int_{\mathbb{R}}\int_0^l|\Delta
\bar{v}|^2\,dxdy)^{\frac{1}{2}},\\
\label{3.2.18}\int_{\mathbb{R}}\int_0^l|\bar{v}||\Delta \bar{v}|\,dxdy
\leq &C(\int_{\mathbb{R}}\int_0^l|
\bar{v}|^2\,dxdy)^{\frac{1}{2}}(\int_{\mathbb{R}}\int_0^l|\Delta
\bar{v}|^2\,dxdy)^{\frac{1}{2}},\\
\label{3.2.19}\int_{\mathbb{R}}\int_0^l(\int_{-h}^{0}|\nabla T|\,dz)|\Delta\bar{v}|\,dxdy\leq &C(\int_{\Omega}| \nabla T|^2\,dxdydz)^{\frac{1}{2}}(\int_{\mathbb{R}}\int_0^l|\Delta
\bar{v}|^2\,dxdy)^{\frac{1}{2}}.
\end{align}
It follows from \eqref{3.2.15}-\eqref{3.2.19} that
\begin{align*}
&\frac{d}{dt}\int_{\mathbb{R}}\int_0^l|
\nabla\bar{v}(x,y,t)|^2\,dxdy+\frac{1}{Re_1}\int_{\mathbb{R}}\int_0^l|\Delta
\bar{v}(x,y,t)|^2\,dxdy\\
\leq &C\int_{\mathbb{R}}\int_0^l| \bar{v}(x,y,t)|^2(\int_{\mathbb{R}}\int_0^l|
\nabla\bar{v}(x,y,t)|^2\,dxdy)^2+C\| \nabla T\|_2^2+C\|
v\|^2_2+C\int_{\Omega}|\nabla \tilde{v}|^2|\tilde{v}|^4\,dxdydz+C\|\nabla\tilde{v}\|_2^2.
\end{align*}

   In view of \eqref{3.2.3}-\eqref{3.2.4}, \eqref{3.2.14} and the uniform Gronwall inequality, we obtain
\begin{align}\label{3.2.20}
\int_{\mathbb{R}}\int_0^l|\nabla\bar{v}(x,y,t)|^2\,dxdy\leq\rho_5
\end{align}
for any $t\geq T_2+3.$
\subsubsection{$H_1$ estimates of $v_{z}$}
 Denote $u=v_{z}$. It is clear that $u$ satisfies the following
 equation obtained by differentiating the first equation \eqref{2.4} with respect to
 $z$:
 \begin{align}\label{3.2.21}
\frac{\partial u}{\partial t}+L_1u+(v\cdot\nabla)u-(\int_{-h}^z\nabla\cdot v(x,y,r,t)dr)\frac{\partial u}{\partial
 z}+(u\cdot\nabla)v-(\nabla\cdot v)u+\frac{1}{Ro}f\vec{k}\times
 u-\nabla T=0
 \end{align}
 with the boundary conditions
 \begin{align}\label{3.2.22}
u|_{\Gamma_u}=0,u|_{\Gamma_b}=0,u(x,0,z,t)=u(x,l,z,t)=0.
\end{align}
Multiplying \eqref{3.2.21} by $u$ and integrating over $\Omega,$ we obtain
\begin{align*}
\frac{1}{2}\frac{d}{dt}\|
u(t)\|^2_2+\|u(t)\|^2=&-\int_{\Omega}[(u\cdot\nabla)v-(\nabla\cdot v)u -\nabla
T]\cdot u\,dxdydz\\
\leq& C\int_{\Omega}|v||u||\nabla u|\,dxdydz +C\int_{\Omega}|T||\nabla u|\,dxdydz\\
\leq&\|T\|_2\|\nabla u\|_2+\|v\|_6\| u\|_3\|\nabla u\|_2,
\end{align*}
which entails
\begin{align*}
\frac{d}{dt}\|u(t)\|^2_2+\|u(t)\|^2\leq C\|v\|_6^4\|u\|^2_2+C\|T\|^2_2.
\end{align*}

It is shown in \cite{jn} that
\begin{align*}
\|v\|_6\leq Ch^{-\frac{1}{3}}\|v\|_2+Ch^{\frac{1}{6}}\left(\int_{\mathbb{R}}\int_0^l|\nabla\bar{v}(x,y)|^2\,dxdy\right)^{\frac{1}{2}}+\|\tilde{v}\|_6,
\end{align*}
which implies that
\begin{align}\label{3.2.23}
\|v(t)\|_6\leq \rho_6
\end{align}
for any $t\geq T_{2}+3.$

From the uniform Gronwall inequality, \eqref{3.2.3}-\eqref{3.2.4} and
\eqref{3.2.23}, we deduce
\begin{align}\label{3.2.24}
\|\partial_zv(t)\|^2_2+\int_t^{t+1}\|\partial_zv(r)\|^2\,dr\leq \rho_7
\end{align}
for any $t\geq T_2+4$.
\subsubsection{$(L^2(\Omega))^2$ estimates of $\nabla v$}
Taking the $L^2(\Omega)$ inner product of the first equation of \eqref{2.4} with $-\Delta v$ and using \eqref{2.5}, we deduce
\begin{align}\label{3.2.25}
\nonumber&\frac{1}{2}\frac{d}{dt}\|\nabla v(t)\|^2_2+\frac{1}{Re_1}\|\Delta
v\|_2^2+\frac{1}{Re_2}\|\nabla v_z\|_2^2\\
\nonumber\leq
&C\int_{\mathbb{R}}\int_0^l(\int_{-h}^0|\nabla v|\,dz)(\int_{-h}^0|\partial_z
v||\Delta v|\,dz)\,dxdy+C\int_{\mathbb{R}}\int_0^l(\int_{-h}^0|\nabla
T|\,dz)(\int_{-h}^0|\Delta v|\,dz)\,dxdy\\
\nonumber&+C\int_{\Omega}|v||\nabla v||\Delta v|\,dxdydz+C\int_{\Omega}|v
||\Delta v|\,dxdydz\\
\leq
&C\|v_z\|_2^{\frac{1}{2}}\|\nabla v_z\|_2^{\frac{1}{2}}\|\nabla v\|_2^{\frac{1}{2}}\|\Delta v\|_2^{\frac{3}{2}}+C\|\nabla
T\|_2\|\Delta v\|_2
+C\|v\|_6\|\nabla
v\|^{\frac{1}{2}}_2\|\Delta
v\|_2^{\frac{3}{2}}+C\|v\|_2\|\Delta v\|_2,
\end{align}
which implies that
\begin{align*}
\frac{d}{dt}\|\nabla
v(t)\|^2_2+\frac{1}{Re_1}\|\Delta
v\|_2^2+\frac{1}{Re_2}\|\nabla v_z(t)\|_2^2
\leq
C(\|v\|^4_6+\|\partial_zv\|^2_2\|\nabla\partial_zv\|^2_2)\|\nabla v\|^2_2+C\|\nabla T\|^2_2+C\| v\|^2_2.
\end{align*}

 From \eqref{3.2.3}-\eqref{3.2.4}, \eqref{3.2.23}-\eqref{3.2.24} and the uniform Gronwall inequality, we deduce
\begin{align}\label{3.2.26}
\|\nabla v(t)\|^2_2+\frac{1}{Re_1}\int_t^{t+1}\|\Delta
v(r)\|_2^2\,dr+\frac{1}{Re_2}\int_t^{t+1}\|\nabla
v_z(r)\|_2^2\,dr\leq\rho_8
\end{align}
for any $t\geq T_2+5.$
\subsubsection{$V_2$ estimates of $T$}
Taking the $L^2(\Omega)$ inner product of the second equation of \eqref{2.4} with $L_2T$ and combining \eqref{2.5}, we find
\begin{align*}
&\frac{1}{2}\frac{d}{dt}\|T(t)\|^2+\|L_2T(t)\|_2^2\\
\leq &C\int_{\Omega}|v||\nabla T||L_2T|\,dxdydz+C\int_{\Omega}|Q ||L_2T|\,dxdydz+C\int_{\mathbb{R}}\int_0^l(\int_{-h}^0|\nabla v|\,dz)(\int_{-h}^{0}|\partial_z
T||L_2T|\,dz)\,dxdy\\
\leq &C\|v\|_6\|\nabla T\|_3\|L_2T\|_2+C\|Q\|_2\|L_2T\|_2+C\|\nabla
v\|_2^{\frac{1}{2}}\|\Delta
v\|_2^{\frac{1}{2}}\|\partial_z
T\|_2^{\frac{1}{2}}\|\nabla\partial_z T\|_2^{\frac{1}{2}}\|L_2T\|_2,
\end{align*}
which entails
\begin{align}\label{3.2.27}
\frac{1}{2}\frac{d}{dt}\|T(t)\|^2+\|L_2T(t)\|_2^2\leq C\|v\|_6^4\|\nabla T\|^2_2+C\|\nabla
v\|^2_2\|\Delta v\|^2_2\|\partial_zT\|^2_2+C\|Q\|^2_2.
\end{align}
We deduce from \eqref{3.2.3}-\eqref{3.2.4}, \eqref{3.2.23}-\eqref{3.2.24}, \eqref{3.2.26} and the uniform Gronwall inequality that
\begin{align}\label{3.2.28}
\|T(t)\|^2+\int_t^{t+1}\|L_2T(r)\|^2_2\,dr\leq\rho_9
\end{align}
for any $t\geq T_2+6.$
\subsubsection{$(L^6(\Omega))^2$ estimates of $\partial_z v$}
Taking the $L^{2}(\Omega)$ inner product of equation
\eqref{3.2.21} with $|u|^4u$ and combining \eqref{3.2.22}, we obtain
\begin{align}\label{3.2.29}
\nonumber&\frac{1}{6}\frac{d}{dt}\|u(t)\|_6^6+\frac{1}{Re_1}\||\nabla u||u|^2\|_2^2+\frac{1}{Re_2}\||\partial_z u||u|^2\|_2^2+\frac{4}{9}\||u|^3\|^2\\
=&\int_{\Omega}(\nabla\cdot v)|u|^6\,dxdydz+\int_{\Omega}\nabla T\cdot |u|^4u\,dxdydz-\int_{\Omega}[(u\cdot\nabla)v]\cdot |u|^4u\,dxdydz.
\end{align}
H\"{o}lder inequality and Sobolev embedding theorem entail
\begin{align}
\nonumber\label{3.2.30}\left|\int_{\Omega}(\nabla\cdot v)|u|^6\,dxdydz\right|\leq &C\int_{\Omega}|v||\nabla |u|^3||u|^3\,dxdydz\\
\nonumber\leq&C\|v\|_6\||u|^3\|_3\|\nabla |u|^3\|_2\\
\leq&C\|v\|_6\|u\|_6^{\frac{3}{2}}(\|\nabla|u|^3\|_2+\|\partial_z|u|^3\|_2)^{\frac{3}{2}},\\
\nonumber\label{3.2.31}\left|\int_{\Omega}\nabla T\cdot |u|^4u\,dxdydz\right|\leq &C\int_{\Omega}|T||\nabla |u|^3||u|^2\,dxdydz+C\int_{\Omega}|T||\nabla u||u|^4\,dxdydz\\
\leq&C\|T\|_6\|\nabla |u|^3\|_2\|u\|_6^2+C\|T\|_6\||\nabla u||u|^2\|_2\|u\|_6^2
\end{align}
and
\begin{align}\label{3.2.32}
\nonumber\left|-\int_{\Omega}[(u\cdot\nabla)v]\cdot |u|^4u\,dxdydz\right|\leq &C\int_{\Omega}|u|^5|v||\nabla u|\,dxdydz+C\int_{\Omega}|u|^3|v||\nabla |u|^3|\,dxdydz\\
\leq&C\|v\|_6\||u|^3\|_3\|\nabla|u|^3\|_2+C\|v\|_6\||u|^3\|_3\||\nabla u||u|^2\|_2.
\end{align}
We deduce from \eqref{3.2.29}-\eqref{3.2.32} that
\begin{align}\label{3.2.33}
\frac{d}{dt}\|u(t)\|_6^2\leq C\|v\|_6^4\|u\|_6^2+C\|T\|_6^2.
\end{align}

Thanks to \eqref{3.2.3}, \eqref{3.2.23}-\eqref{3.2.24} and the uniform Gronwall inequality, we find
\begin{align}\label{3.2.34}
\|\partial_z v\|_6^2\leq \rho_{10}
\end{align}
for any $t\geq T_2+7.$
\subsubsection{$L^6(\Omega)$ estimates of $\partial_z T$}
Denote $\theta=T_z$. It is clear that $\theta$ satisfies the following
 equation obtained by differentiating the second equation of \eqref{2.4} with respect to
 $z$:
 \begin{align}\label{3.2.35}
\frac{\partial \theta}{\partial t}+L_2\theta+v\cdot\nabla\theta-\left(\int_{-h}^z\nabla\cdot v(x,y,r,t)dr\right)\frac{\partial \theta}{\partial
 z}+\partial_zv\cdot\nabla T-(\nabla\cdot v)\theta =Q_z
 \end{align}
 with the boundary conditions
 \begin{align}\label{3.2.36}
(\frac{1}{Rt_2}\theta+\alpha T)|_{\Gamma_u}=0,\theta|_{\Gamma_b}=0,\frac{\partial \theta}{\partial\vec{n}}|_{\Gamma_l}=0.
\end{align}
Taking the $L^2(\Omega)$ inner product of equation
\eqref{3.2.35} with $|\theta|^4\theta,$ we obtain
\begin{align}\label{3.2.37}
\nonumber&\frac{1}{6}\frac{d}{dt}\|\theta(t)\|_6^6+\frac{5}{9Rt_1}\int_{\Omega}|\nabla|\theta|^3|^2\,dxdydz+\frac{5}{9Rt_2}\int_{\Omega}|\partial_z|\theta|^3|^2\,dxdydz+\alpha^5 Rt_2^4\int_{\Gamma_u}\frac{\partial \theta}{\partial z}|T|^4T\,dxdy\\
=&\int_{\Omega}(\nabla\cdot v)|\theta|^6\,dxdydz+\int_{\Omega}Q_z|\theta|^4\theta\,dxdydz-\int_{\Omega}(v_z\cdot\nabla T)|\theta|^4\theta\,dxdydz.
 \end{align}

In the following, we estimate each term in the right hand side of \eqref{3.2.37}.
\begin{align}
\nonumber\label{3.2.38}\left|\int_{\Omega}(\nabla\cdot v)|\theta|^6\,dxdydz\right|\leq &C\int_{\Omega}|v||\nabla |\theta|^3||\theta|^3\,dxdydz\\
\nonumber\leq&C\|v\|_6\||\theta|^3\|_3\|\nabla |\theta|^3\|_2\\
\leq&C\|v\|_6\||\theta|^3\|_2^{\frac{1}{2}}\||\theta|^3\|^{\frac{3}{2}},\\
\nonumber\label{3.2.39}\left|\int_{\Omega}Q_z|\theta|^4\theta\,dxdydz\right|\leq &C\int_{\Omega}|Q_z||\theta|^5\,dxdydz\\
\nonumber\leq&C\|Q_z\|_2\||\theta|^3\|_{\frac{10}{3}}^{\frac{5}{3}}\\
\leq&C\|Q\|\||\theta|^3\|_2^{\frac{2}{3}}\||\theta|^3\|,\\
\nonumber\label{3.2.40}\left|-\int_{\Omega}(v_z\cdot\nabla T) |\theta|^4\theta\,dxdydz\right|\leq &C\int_{\Omega}|v_z||\nabla T| |\theta|^5\,dxdydz\\
\nonumber\leq&C\|v_z\|_6\|\nabla T\|_3\||\theta|^3\|_{\frac{10}{3}}^{\frac{5}{3}}\\
\leq&C\|v_z\|_6\|\nabla T\|_2^{\frac{1}{2}}\|\nabla T\|^{\frac{1}{2}}\||\theta|^3\|_2^{\frac{2}{3}}\||\theta|^3\|,\\
\nonumber\label{3.2.41}\alpha^5 Rt_2^4\int_{\Gamma_u}\frac{\partial \theta}{\partial z}|T|^4T\,dxdy=&\alpha^5 Rt_2^5\int_{\Gamma_u}(\frac{\partial T}{\partial t}+v\cdot\nabla T-\frac{1}{Rt_1}\Delta T-Q)|T|^4T\,dxdy\\
\nonumber=&\frac{\alpha^5 Rt_2^5}{6}\frac{d}{dt}\int_{\Gamma_u}|T|^6\,dxdy+\alpha^5 Rt_2^5\int_{\Gamma_u}(v\cdot\nabla T)|T|^4T\,dxdy\\
&+\frac{5\alpha^5 Rt_2^5}{9Rt_1}\int_{\Gamma_u}|\nabla |T|^3|^2\,dxdy-\alpha^5 Rt_2^5\int_{\Gamma_u}Q|T|^4T\,dxdy,\\
\nonumber\label{3.2.42}&\left|\alpha^5 Rt_2^5\int_{\Gamma_u}(v\cdot\nabla T)|T|^4T\,dxdy-\alpha^5 Rt_2^5\int_{\Gamma_u}Q|T|^4T\,dxdy\right|\\
\nonumber\leq&C\|v\|_{L^4(\Gamma_u)}\|\nabla |T|^3\|_{L^2(\Gamma_u)}\||T|^3\|_{L^4(\Gamma_u)}+C\|Q\|_{L^2(\Gamma_u)}\||T|^3\|_{L^{\frac{10}{3}}(\Gamma_u)}^{\frac{5}{3}}\\
\leq&C\|v\|\|\nabla |T|^3\|_{L^2(\Gamma_u)}^{\frac{3}{2}}\|T\|_{L^6(\Gamma_u)}^{\frac{3}{2}}+C\|Q\|\||T|^3\|_{L^2(\Gamma_u)}\||T|^3\|_{H^1(\Gamma_u)}^{\frac{2}{3}}.
\end{align}
It follows from \eqref{3.2.37}-\eqref{3.2.42} that
\begin{align}\label{3.2.43}
\nonumber&\frac{d}{dt}(\|\theta(t)\|_6^6+\alpha^5 Rt_2^5\|T(t)\|_{L^6(\Gamma_u)}^6)+\frac{2}{Rt_1}\int_{\Omega}|\nabla|\theta|^3|^2\,dxdydz+\frac{2}{Rt_2}\int_{\Omega}|\partial_z|\theta|^3|^2\,dxdydz+\frac{2\alpha^5 Rt_2^5}{Rt_1}\int_{\Gamma_u}|\nabla |T|^3|^2\,dxdy\\
\leq&C\|v\|_6^4\|\theta\|_6^6+C\|Q\|^2\|\theta\|_6^4+C\|v_z\|_6^2\|\nabla T\|_2\|\nabla T\|\|\theta\|_6^4+C\|v\|^4\|T\|_{L^6(\Gamma_u)}^6+C\|Q\|^{\frac{3}{2}}\|T\|_{L^6(\Gamma_u)}^{\frac{9}{2}}.
\end{align}

 As a result of \eqref{3.2.23}-\eqref{3.2.24}, \eqref{3.2.26}, \eqref{3.2.28}, \eqref{3.2.34} and the uniform Gronwall inequality, we find
\begin{align}\label{3.2.44}
\|\partial_z T(t)\|_6^2+\alpha Rt_2\|T(t)\|_{L^6(\Gamma_u)}^2\leq \rho_{11}
\end{align}
for any $t\geq T_2+8.$
\subsubsection{$H$ estimates of $(\partial_t v,\partial_t T)$ }
 Denote $\pi=\partial_t v,$ $\xi=\partial_t T.$ It is clear that $\pi,$ $\xi$ satisfy the following
 equations obtained by differentiating the first equation and the second equation of \eqref{2.4} with respect to
 $t,$ respectively:
\begin{equation}\label{3.2.45}
\begin{cases}
&\frac{\partial \pi}{\partial t}+L_1\pi+(v\cdot\nabla)\pi-\left(\int_{-h}^z\nabla\cdot
v(x,y,r,t)dr\right)\frac{\partial \pi}{\partial
z}+\frac{1}{Ro}f\vec{k}\times \pi+\nabla\partial_t
p_s(x,y,t)\\
& +(\pi\cdot\nabla)v-\left(\int_{-h}^z\nabla\cdot
\pi(x,y,r,t)dr\right)\frac{\partial v}{\partial z}
-\int_0^z \nabla\xi(x,y,r,t)dr=0,\\
&\frac{\partial \xi}{\partial
t}+L_2\xi+v\cdot\nabla\xi-\left(\int_{-h}^z\nabla\cdot
v(x,y,r,t)dr\right)\frac{\partial \xi}{\partial
z}+\pi\cdot\nabla T-\left(\int_{-h}^z\nabla\cdot \pi(x,y,r,t)dr\right)\frac{\partial
T}{\partial z}=0
\end{cases}
\end{equation}
with the boundary conditions
\begin{equation}\label{3.2.46}
\begin{cases}
&\frac{\partial \pi}{\partial z}|_{\Gamma_u}=0,\frac{\partial\pi}{\partial z}|_{\Gamma_b}=0,\pi\cdot\vec{n}|_{\Gamma_l}=0,\frac{\partial \pi}{\partial\vec{n}}\times\vec{n}|_{\Gamma_l}=0,\\
&(\frac{1}{Rt_2}\frac{\partial \xi}{\partial z}+\alpha \xi)|_{\Gamma_u}=0,\frac{\partial \xi}{\partial z}|_{\Gamma_b}=0,\frac{\partial \xi}{\partial\vec{n}}|_{\Gamma_l}=0.
\end{cases}
\end{equation}
Multiplying the first equation and the second equation of \eqref{3.2.45} by $\pi,$ $\xi,$ respectively, and integrating over $\Omega,$ we obtain
\begin{align}\label{3.2.47}
\nonumber&\frac{1}{2}\frac{d}{dt}\|\pi(t)\|_2^2+\|\pi(t)\|^2\\
\nonumber=& -\int_{\Omega}\left[(\pi\cdot\nabla)v-\left(\int_{-h}^z\nabla\cdot
\pi(x,y,r,t)dr\right)\frac{\partial v}{\partial z}\right]\cdot\pi\,dxdydz+\int_{\Omega}\int_0^z \nabla\xi(x,y,r,t)dr\cdot\pi\,dxdydz\\
\leq&\|v\|_6\|\pi\|_3\|\nabla\pi\|_2+C\|\nabla\pi\|_2^{\frac{3}{2}}\|\partial_z v\|_2^{\frac{1}{2}}\|\nabla\partial_z v\|_2^{\frac{1}{2}}\|\pi\|_2^{\frac{1}{2}}+C\|\xi\|_2\|\nabla\pi\|_2
\end{align}
and
\begin{align}\label{3.2.48}
\nonumber\frac{1}{2}\frac{d}{dt}\|\xi(t)\|_2^2+\|\xi(t)\|^2
=&-\int_{\Omega}(\pi\cdot\nabla T)\xi\,dxdydz
+\int_{\Omega}\left(\int_{-h}^z\nabla\cdot \pi(x,y,r,t)dr\right)\frac{\partial
T}{\partial z}\xi\,dxdydz\\
\leq&C\|\nabla\pi\|_2\|\partial_z T\|_2^{\frac{1}{2}}\|\nabla\partial_z T\|_2^{\frac{1}{2}}\|\xi\|_2^{\frac{1}{2}}\|\nabla\xi\|_2^{\frac{1}{2}}
+C\|\pi\|_3\|T\|_6\|\nabla\xi\|_2+C\|\nabla\pi\|_2\|T\|_6\|\xi\|_3.
\end{align}
It follows from \eqref{3.2.47}-\eqref{3.2.48} and Young inequality that
\begin{align}\label{3.2.49}
\frac{d}{dt}(\|\xi(t)\|_2^2+\|\pi(t)\|_2^2)\leq A(t)(\|\xi\|_2^2+\|\pi\|_2^2),
\end{align}
where
\begin{align*}
A(t)=&C(1+\|T\|_6^4+\|v\|_6^4+\|\partial_z
T\|_2^2\|\nabla
\partial_zT\|_2^2+\|\partial_zv\|_2^2\|\nabla
\partial_zv\|_2^2).
\end{align*}
Multiplying the second equation of \eqref{2.4} by $\xi$ and integrating over $\Omega,$ we obtain
\begin{align}\label{3.2.50}
\nonumber\|\xi\|_2^2=&-\int_{\Omega}\left[v\cdot\nabla T-\left(\int_{-h}^z
\nabla\cdot v(x,y,r,t)\,dr\right)\frac{\partial T}{\partial
z}\right]\xi
\,dxdydz-\int_{\Omega}L_2 T\xi\,dxdydz+\int_{\Omega}Q\xi\, dxdydz\\
\leq&\|Q\|_2\|\xi\|_2+C\|v\|_6\|\nabla
T\|_3\|\xi\|_2+\|L_2T\|_2\|\xi\|_2+C\|\nabla
v\|_2^{\frac{1}{2}}\|\Delta v\|_2^{\frac{1}{2}}\|\partial_z T\|_2^{\frac{1}{2}}\|\nabla T_z\|_2^{\frac{1}{2}}\|\xi\|_2.
\end{align}
Similarly, we have
\begin{align}\label{3.2.51}
\nonumber\|\pi\|_2^2=&-\int_{\Omega}\left[(v\cdot\nabla)v-\left(\int_{-h}^z\nabla\cdot v(x,y,r,t)\,dr\right)\frac{\partial v}{\partial z}\right]\cdot\pi\,dxdydz\\
\nonumber&-\int_{\Omega}\left(\frac{f}{Ro}\vec{k}\times v+\int_0^z\nabla T(x,y,r,t)\,dr+L_1v\right)\cdot\pi\,dxdydz\\
\leq&C\|v\|_6\|\nabla v\|_3\|\pi\|_2+C\|\nabla v\|_2^{\frac{1}{2}}\|\Delta v\|_2^{\frac{1}{2}}\|v_z\|_2^{\frac{1}{2}}\|\nabla v_z\|_2^{\frac{1}{2}}\|\pi\|_2+C\|v\|_2\|\pi\|_2+C\|\nabla T\|_2\|\pi\|_2+\|L_1v\|_2\|\pi\|_2.
\end{align}
We derive from \eqref{3.2.50}-\eqref{3.2.51} and Young inequality that
\begin{align}\label{3.2.52}
\nonumber\|\pi\|_2^2+\|\xi\|_2^2\leq&C(1+\|v\|_6^4)(\|\nabla
T\|_2^2+\|\nabla v\|_2^2)+C\|Q\|_2^2+C(1+\|\partial_z
T\|_2^2)\|L_2T\|_2^2\\
&+C(1+\|\nabla v\|_2^2+\|\partial_z v\|_2^2)\|L_1 v\|_2^2.
\end{align}
By use of \eqref{3.2.23}-\eqref{3.2.24}, \eqref{3.2.26}, \eqref{3.2.28}, \eqref{3.2.34}, \eqref{3.2.44} and the uniform Gronwall inequality, we obtain
\begin{align}\label{3.2.53}
\|v_t(t)\|_2^2+\|T_t(t)\|_2^2+\int_t^{t+1}\|v_t(r)\|^2+\|T_t(r)\|^2\,dr\leq \rho_{12}
\end{align}
for any $t\geq T_2+9.$
\subsubsection{$H$ estimates of $(L_1v, L_2T)$}
Multiplying the second equation of \eqref{2.4} by $L_2T$ and integrating over $\Omega,$ we obtain
\begin{align}\label{3.2.54}
\nonumber\|L_2T\|_2^2=&-\int_{\Omega}(v\cdot\nabla T)L_2T\,dxdydz+\int_{\Omega}(\int_{-h}^z
\nabla\cdot v(x,y,r,t)\,dr)\frac{\partial T}{\partial z}L_2T\,dxdydz\\
\nonumber&-\int_{\Omega}\partial_t TL_2T\,dxdydz+\int_{\Omega}QL_2T\, dxdydz\\
\nonumber\leq&C\|\nabla v\|_2^{\frac{1}{2}}\|\Delta v\|_2^{\frac{1}{2}}\|T_z\|_6\|L_2T\|_2+C\|\nabla v\|_2\|T_z\|_6\|L_2T\|_2\\
&+\|Q\|_2\|L_2T\|_2+C\|v\|_6\|\nabla
T\|_3\|L_2T\|_2+\|\partial_t T\|_2\|L_2T\|_2,
\end{align}
Similarly, we have
\begin{align}\label{3.2.55}
\nonumber&\|L_1v\|_2^2=-\int_{\Omega}\left[(v\cdot\nabla)v-\left(\int_{-h}^z\nabla\cdot v(x,y,r,t)\,dr\right)\frac{\partial v}{\partial z}\right]\cdot L_1v\,dxdydz\\
\nonumber&-\int_{\Omega}\left(\frac{f}{Ro}\vec{k}\times v-\int_0^z\nabla T(x,y,r,t)\,dr+\partial_t v\right)\cdot L_1v\,dxdydz\\
\nonumber\leq&C\|v\|_6\|\nabla v\|_3\|L_1v\|_2+C\|\nabla v\|_2^{\frac{1}{2}}\|\Delta v\|_2^{\frac{1}{2}}\|v_z\|_6\|L_1v\|_2+C\|\nabla v\|_2\|v_z\|_6\|L_1v\|_2\\
&+C\|v\|_2\|L_1v\|_2+C\|\nabla T\|_2\|L_1v\|_2+\|v_t\|_2\|L_1v\|_2.
\end{align}
It follows from \eqref{3.2.54}-\eqref{3.2.55} and Young inequality that
\begin{align}\label{3.2.56}
\nonumber&\|L_1v\|_2^2+\|L_2T\|_2^2\\
\leq&C(1+\|v\|_6^4+\|T_z\|_6^4+\|v_z\|_6^4)(\|\nabla
T\|_2^2+\|\nabla v\|_2^2)+C(\|v_t\|_2^2+\|T_t\|_2^2)+C\|Q\|_2^2.
\end{align}
By virtue of \eqref{3.2.23}-\eqref{3.2.24}, \eqref{3.2.26}, \eqref{3.2.28}, \eqref{3.2.34}, \eqref{3.2.44}, \eqref{3.2.53}, we conclude
\begin{align}\label{3.2.57}
\|L_1v(t)\|_2^2+\|L_2T(t)\|_2^2\leq\rho_{13}
\end{align}
for any $t\geq T_2+9.$ %Moreover, we obtain
%\begin{align}\label{3.2.58}
%\int_{\mathbb{R}}\int_0^l|p_s(x,y,t)|^2\,dxdy\leq\rho_{14}
%\end{align}
%for any $t\geq T_2+9.$
\section{The existence of a global attractor}
\def\theequation{4.\arabic{equation}}\makeatother
\setcounter{equation}{0}
\indent In this section, we will prove the existence of a global
attractor in $V$ for problem \eqref{2.4}-\eqref{2.6}.

First of all, from \eqref{3.2.3}-\eqref{3.2.4}, \eqref{3.2.24},\eqref{3.2.26}, \eqref{3.2.28}, \eqref{3.2.57}, we immediately conclude the
following result.
\begin{theorem}\label{4.4}
Assume that $Q \in H^1(\Omega).$ Then the semigroup $\{S(t)\}_{t \geq
0}$ generated by problem \eqref{2.4}-\eqref{2.6} possesses an absorbing set in $(H^2(\Omega))^3\cap V.$
\end{theorem}
Next, we prove the asymptotical compactness of the semigroup $\{S(t)\}_{t \geq 0}$
generated by problem \eqref{2.4}-\eqref{2.6}.
\begin{theorem} \label{4.5}
Assume that $Q \in H^1(\Omega).$ Let $\{S(t)\}_{t \geq 0}$ be the semigroup
generated by problem \eqref{2.4}-\eqref{2.6} and let $\{(v_{0n},T_{0n})\}_{n=1}^{\infty}$ be a sequence in $V$ which converges weakly in $V$ to an element $(v_0,T_0)\in V,$ denote by $(v^n(t),T^n(t))=S(t)(v_{0n},T_{0n})$ and $(v(t),T(t))=S(t)(v_0,T_0)$ for any $t\geq 0.$ Then
\begin{align}
v^n(t)\rightharpoonup v(t)\,\,\,\textit{weakly in}\,\,H_1,\,\,\,\forall\,\,\,t\geq 0
\end{align}
and
\begin{align}
(v^n(t),T^n(t))\rightharpoonup (v(t),T(t))\,\,\,\textit{weakly in}\,\,L^2(0,\mathcal{T};V),\,\,\,\forall\,\,\,\mathcal{T}> 0.
\end{align}
\end{theorem}
\textbf{Proof.} It is easily verified that
\begin{align}\label{4.6}
(v^n(t),T^n(t))\,\,\,\textit{is bounded in}\,\,L^{\infty}(\mathbb{R}^+;V)\cap L^2(0,\mathcal{T};V\cap (H^2(\Omega))^3),\,\,\,\forall\,\,\,\mathcal{T}> 0
\end{align}
and
\begin{align}\label{4.7}
(v^n_t(t),T^n_t(t))\,\,\,\textit{is bounded in}\,\,L^2(0,\mathcal{T};H),\,\,\,\forall\,\,\,\mathcal{T}> 0.
\end{align}
Then, for any $u\in H_1$ and $0\leq t\leq t+a\leq\mathcal{T}$ with $\mathcal{T}>0,$
\begin{align}\label{4.8}
\nonumber\int_\Omega(v^n(t+a)-v^n(t))\cdot u\,dxdydz=&\int_\Omega\int_t^{t+a}v^n_t(s)\cdot u\,dxdydzds\\
\nonumber\leq&\int_t^{t+a}\|v^n_t(s)\|_{H_1}\|u\|_2\,ds\\
\nonumber\leq&\|v^n_t\|_{L^2(0,\mathcal{T};H_1)}\|u\|_2a^{\frac{1}{2}}\\
\leq&C_\mathcal{T}a^{\frac{1}{2}}\|u\|_2,
\end{align}
where $C_\mathcal{T}$ is a positive constant independent of $n.$ Substituting $u=v^n(t+a)-v^n(t)$ into equation \eqref{4.8}, we obtain
\begin{align*}
\|v^n(t+a)-v^n(t)\|_2^2
\leq C_\mathcal{T}a^{\frac{1}{2}}\|v^n(t+a)-v^n(t)\|_2,
\end{align*}
which entails
\begin{align*}
\int_0^{\mathcal{T}-a}\|v^n(t+a)-v^n(t)\|_2^2\,dt\leq\mathcal{T}C^2_\mathcal{T}a
\end{align*}
for another positive constant independent of $n.$ Therefore, we obtain
\begin{align*}
\lim_{a\rightarrow 0}\sup_n\int_0^{\mathcal{T}-a}\|v^n(t+a)-v^n(t)\|_{L^2(\Omega_r)}^2\,dt=0
\end{align*}
for any $r>0,$ where $\Omega_r=\{(x,y,z)\in\Omega: |x|\leq r\}.$ From \eqref{4.6}-\eqref{4.7} and the Aubin-Lions Compactness Lemma, we deduce that
\begin{equation}\label{4.9}
\begin{cases}
&(v^n|_{\Omega_r}(t),T^n|_{\Omega_r}(t))\,\,\,\textit{is relatively compact in}\,\,L^2(0,\mathcal{T};(L^2(\Omega_r))^3),\,\,\,\forall\,\,\,\mathcal{T}> 0,\,\,\,\forall\,\,\,r> 0,\\
&(v^n|_{\Omega_r}(t),T^n|_{\Omega_r}(t))\,\,\,\textit{is relatively compact in}\,\,L^2(0,\mathcal{T};(H^1(\Omega_r))^3),\,\,\,\forall\,\,\,\mathcal{T}> 0,\,\,\,\forall\,\,\,r> 0.
\end{cases}
\end{equation}
We infer from \eqref{4.6}-\eqref{4.7}, \eqref{4.9} and the diagonal process that there exists some subsequence $\{(v^{n_j},T^{n_j})\}_{j=1}^{\infty}$ of the sequence $\{(v^n,T^n)\}_{n=1}^{\infty}$ and some $(v_1,T_1)\in L^{\infty}(\mathbb{R}^+;V)\cap L^2_{loc}(\mathbb{R}^+;(H^2(\Omega))^3\cap V)$ such that
\begin{equation}\label{4.10}
\begin{cases}
&(v^{n_j},T^{n_j})\rightarrow (v_1,T_1)\,\,\,\textit{weakly-star in}\,\,L^{\infty}(0,\mathcal{T};V),\\
&(v^{n_j},T^{n_j})\rightarrow (v_1,T_1)\,\,\,\textit{weakly in}\,\,L^2_{loc}(\mathbb{R}^+;(H^2(\Omega))^3\cap V),\\
&(v^{n_j},T^{n_j})\rightarrow (v_1,T_1)\,\,\,\textit{strongly in}\,\,L^2_{loc}(\mathbb{R}^+;(L^2(\Omega_r))^3),\,\,\,\forall\,\,\,r> 0,\\
&(v^{n_j},T^{n_j})\rightarrow (v_1,T_1)\,\,\,\textit{strongly in}\,\,L^2_{loc}(\mathbb{R}^+;(H^1(\Omega_r))^3),\,\,\,\forall\,\,\,r> 0,
\end{cases}
\end{equation}
which implies that $(v_1,T_1)$ is a solution of problem \eqref{2.4}-\eqref{2.6} with $(v_1(0),T_1(0))=(v_0,T_0)\in V.$ By the uniqueness of strong solution of problem problem \eqref{2.4}-\eqref{2.6}, we obtain $v_1=v.$ Then by a contradiction argument, we deduce that the whole sequence $\{v^n\}_{n=1}^{\infty}$ converges to $v$ in the sense of \eqref{4.10}.

It follows from \eqref{4.10} that $\{v^n(t)\}_{n=1}^{\infty}$ converges strongly in $(L^2(\Omega_r))^2$ to $v(t)$ for almost every $t\geq 0$ and any $r>0.$ Then for any $u\in\mathcal{V}_1,$
\begin{align}\label{4.11}
\int_\Omega v^n(t)\cdot u\,dxdydz\rightarrow\int_\Omega v(t)\cdot u\,dxdydz,\,\,a.e.t\in\mathbb{R}^+.
\end{align}
From \eqref{4.6} and \eqref{4.8}, we deduce that $\{\int_\Omega v^n(t)\cdot u\,dxdydz\}_{n=1}^{\infty}$ is equibounded and equicontinuous on $[0,\mathcal{T}]$ for any $\mathcal{T}>0.$ Therefore, we have
\begin{align}\label{4.12}
\int_\Omega v^n(t)\cdot u\,dxdydz\rightarrow\int_\Omega v(t)\cdot u\,dxdydz,\,\,\forall\,t\in\mathbb{R}^+,\,\,\,\forall\,\,u\in\mathcal{V}_1,
\end{align}
which entails
\begin{align*}
v^n(t)\rightharpoonup v(t)\,\,\,\textit{weakly in}\,\,H,\,\,\,\forall\,\,\,t\geq 0.
\end{align*}
\qed\hfill
\begin{theorem} \label{4.13}
Assume that $Q \in H^1(\Omega).$ Then the semigroup $\{S(t)\}_{t \geq 0}$
generated by problem \eqref{2.4}-\eqref{2.6} is asymptotically compact in $V.$
\end{theorem}
\textbf{Proof.} Choose a smooth function $\eta$ such that $0\leq
\eta(s)\leq 1$ for all $s\in \mathbb{R}^{+},$ and
\begin{equation*}
\eta(s)=
\begin{cases}
&0,\,\,\,\, 0\leq s\leq 1;\\
&1,\,\,\,\,\, s\geq 2.
\end{cases}
\end{equation*}
Then there exists a positive constant $C$ such that $|\eta'(s)|\leq C$
for any $s\in \mathbb{R}^+.$ Define a projection operator $P_i:\mathbb{R}^3\rightarrow\mathbb{R}^{3-i}$ by $P_1(u_1,u_2,u_3)=(u_1,u_2)$ and $P_1(u_1,u_2,u_3)=u_3$ for any $(u_1,u_2,u_3)\in\mathbb{R}^3$ and $i=1,2,$ let $S_1(t)=P_1S(t)$ and $S_2(t)=P_2S(t)$ for any $t\geq 0.$

Multiplying the second equation of \eqref{2.4} by $\eta^2(\frac{x^2}{r^2})T$ and integrating over $\Omega,$ we conclude
\begin{align*}
&\frac{1}{2}\frac{d}{dt}\int_\Omega\eta^2(\frac{x^2}{r^2})|T(x,y,z,t)|^2\,dxdydz+\|\eta T(t)\|^2\\
\leq&\int_\Omega\frac{4x^2}{Rt_1r^4}|\eta'(\frac{x^2}{r^2})|^2|T(x,y,z,t)|^2\,dxdydz+\int_\Omega\frac{2|x|}{r^2}|\eta'(\frac{x^2}{r^2})|\eta(\frac{x^2}{r^2})|v(x,y,z,t)||T(x,y,z,t)|^2\,dxdydz\\
&+\int_\Omega\eta^2(\frac{x^2}{r^2})|T(x,y,z,t)||Q(x,y,z)|\,dxdydz\\
\leq&\frac{C}{r^2}\|T\|_2^2+\frac{C}{r}\|v\|_3\|T\|_6\|\eta T\|_2+\|\eta T\|_2\|\eta Q\|_2.
\end{align*}
We infer from Young inequality and Poinc\'{a}re inequality \eqref{3.2.2} that
\begin{align}\label{4.14}
\frac{d}{dt}\int_\Omega\eta^2(\frac{x^2}{r^2})|T(x,y,z,t)|^2\,dxdydz+\frac{1}{2Rt_2h^2+\frac{2h}{\alpha}}\|\eta T\|_2^2
\leq\frac{C}{r^2}(\|T\|_2^2+\|v\|^2\|T\|_6^2)+C\|\eta Q\|_2^2.
\end{align}
From Theorem \ref{4.4}, we conclude that there exists a positive constant $\mathcal{M}$ satisfying for any bounded subset $B\subset V,$ there exists some time $\tau=\tau(B)$ such that
\begin{align}\label{4.15}
\|v(t)\|_{(H^2(\Omega))^2}^2+\|T(t)\|_{H^2(\Omega)}^2\leq\mathcal{M}^2
\end{align}
for any $t\geq\tau.$

It follows from the classical Gronwall inequality, \eqref{4.14}-\eqref{4.15} and $Q\in H^1(\Omega)$ that for any $\epsilon>0,$ there exists some $\mathcal{R}_1$ such that
\begin{align}\label{4.16}
\int_\Omega\eta^2(\frac{x^2}{r^2})|T(x,y,z,t)|^2\,dxdydz\leq C\epsilon
\end{align}
for any $t\geq\tau$ and any $r\geq\mathcal{R}_1.$

Therefore, for any $\epsilon>0,$ we deduce that
\begin{align}\label{4.17}
\int_{\Omega^c_{2\mathcal{R}_1}}|T(x,y,z,t)|^2\,dxdydz\leq C\epsilon
\end{align}
for any $t\geq\tau,$ where $\Omega_{2\mathcal{R}_1}=\{(x,y,z)\in\Omega:|x|\leq 2\mathcal{R}_1\}$ and $E^c=\mathbb{R}^n\backslash E$ for any $E\subset\mathbb{R}^n.$

Let $B$ be any fixed bounded subset in $V,$ for any $\{(v_{0n},T_{0n})\}_{n=1}^{\infty}\subset B$ and any positive sequence $\{t_n\}_{n=1}^{\infty}$ with $t_n\rightarrow+\infty$ as $n\rightarrow+\infty,$ denote by $(v^n(t_n),T^n(t_n))=S(t_n)(v_{0n},T_{0n}),$ we infer from Theorem \ref{4.4} that there exists a bounded ball $B_0$ in $V$ and some time $\tau=\tau(B)>0$ such that
\begin{align*}
S(t)B\subset B_0
\end{align*}
for any $t\geq\tau$ and
\begin{align*}
\|v\|^2+\|T\|^2\leq\rho^2
\end{align*}
for any $(v,T)\in B_0,$ which implies that exists some $K\in\mathbb{N}$ with $t_n\geq\tau$ for any $n\geq K$ such that
\begin{align*}
T^n(t_n)\chi_{\Omega_{2\mathcal{R}_1}}\,\,\,\textit{is uniformly bounded in}\,\,H^1(\Omega_{2\mathcal{R}_1}),
\end{align*}
where $\chi_E$ is an indicator function on $E\subset\mathbb{R}^n.$

Thanks to the compactness of $H^1(\Omega_{2\mathcal{R}_1})\subset L^2(\Omega_{2\mathcal{R}_1}),$ we find that for the above $\epsilon,$ there exists some subsequence $\{T^{n_j}(t_{n_j})\chi_{\Omega_{2\mathcal{R}_1}}\}_{j=1}^{\infty}$ of $\{T^n(t_n)\chi_{\Omega_{2\mathcal{R}_1}}\}_{n=1}^{\infty}$ and some $K_1\in\mathbb{N}$ with $n_j\geq K$ for any $j\geq K_1$ such that
\begin{align*}
\|T^{n_j}(t_{n_j})\chi_{\Omega_{2\mathcal{R}_1}}-T^{n_l}(t_{n_l})\chi_{\Omega_{2\mathcal{R}_1}}\|_{L^2(\Omega_{2\mathcal{R}_1})}\leq\epsilon
\end{align*}
for any $j,l\geq K_1.$ Therefore, we obtain
\begin{align*}
&\|T^{n_j}(t_{n_j})-T^{n_l}(t_{n_l})\|_{L^2(\Omega)}\\
\leq&\|T^{n_j}(t_{n_j})\chi_{\Omega_{2\mathcal{R}_1}}-T^{n_l}(t_{n_l})\chi_{\Omega_{2\mathcal{R}_1}}\|_{L^2(\Omega_{2\mathcal{R}_1})}+\|T^{n_j}(t_{n_j})\chi_{\Omega^c_{2\mathcal{R}_1}}-T^{n_l}(t_{n_l})\chi_{\Omega^c_{2\mathcal{R}_1}})\|_{L^2(\Omega^c_{2\mathcal{R}_1})}\\
\leq&\|T^{n_j}(t_{n_j})\chi_{\Omega_{2\mathcal{R}_1}}-T^{n_l}(t_{n_l})\chi_{\Omega_{2\mathcal{R}_1}}\|_{L^2(\Omega_{2\mathcal{R}_1})}+\|T^{n_j}(t_{n_j})\chi_{\Omega^c_{2\mathcal{R}_1}}\|_{L^2(\Omega^c_{2\mathcal{R}_1})}+\|T^{n_l}(t_{n_l})\chi_{\Omega^c_{2\mathcal{R}_1}}\|_{L^2(\Omega^c_{2\mathcal{R}_1})}\\
\leq&C\epsilon
\end{align*}
for any $j,l\geq K_1.$ Due to $\|T\|\leq C\|T\|_2^\frac{1}{2}\|T\|_{H^2(\Omega)\cap V_2}^\frac{1}{2},$ we know that the sequence $\{T^{n_j}(t_{n_j})\}_{j=1}^{\infty}$ is a Cauchy sequence in $V_2,$ which implies that the semigroup $\{S_2(t)\}_{t \geq 0}$ is asymptotically compact in $V_2.$

Let $\delta=\frac{1}{2lRe_1}$ and define
\begin{align}
a(v,v)=\|v\|^2-\frac{\delta}{2}\|v\|_2^2.
\end{align}
It is easy to prove that
\begin{align}\label{4.18}
\frac{d}{dt}\|v(t)\|^2_2+\delta\|v(t)\|^2_2+2a(v,v)
=2\int_{\Omega}(\int_0^z\nabla T(x,y,r,t)\,dr)\cdot v(x,y,z,t)\,dxdydz
\end{align}
for any solution $(v,T)=(v(t),T(t))=S(t)(v_0,T_0)$ with $(v_0,T_0)\in V.$ From the variation of constants formula, we conclude
\begin{align}\label{4.19}
\nonumber&\|v(t)\|^2_2=e^{-\delta t}\|v_0\|^2_2-2\int_0^te^{-\delta(t-s)}a(v(s),v(s))\,ds\\
+&2\int_0^te^{-\delta(t-s)}\left(\int_{\Omega}(\int_0^z\nabla T(x,y,r,s)\,dr)\cdot v(x,y,z,s)\,dxdydz\right)\,ds.
\end{align}
We are now in position to prove the sequence $\{v^n(t_n)\}_{n=1}^{\infty}$ is a Cauchy sequence in $H_1.$ Thanks to
\begin{align*}
S(t)B\subset B_0
\end{align*}
for any $t\geq\tau,$ we obtain
\begin{align*}
(v^n(t_n),T^n(t_n))\subset B_0
\end{align*}
for any $n\geq K.$ Therefore, the sequence $\{(v^n(t_n),T^n(t_n))\}_{n=1}^{\infty}$ is weakly relatively compact in $V,$ which entails that there exists some subsequence $\{(v^{n_j}(t_{n_j}),T^{n_j}(t_{n_j}))\}_{j=1}^{\infty}$ of $\{(v^n(t_n),T^n(t_n))\}_{n=1}^{\infty}$ and $(v,T)\in B_0$ such that
\begin{align}\label{4.20}
(v^{n_j}(t_{n_j}),T^{n_j}(t_{n_j}))\rightharpoonup (v,T)\,\,\textit{weakly in}\,\,V.
\end{align}
Similarly, for each $J>0,$ we also have
\begin{align}\label{4.21}
(v^n(t_n-J),T^n(t_n-J))\subset B_0
\end{align}
for any $t_n\geq \tau+J.$ Hence, the sequence $\{(v^n(t_n-J),T^n(t_n-J))\}_{n=1}^{\infty}$ is weakly relatively compact in $V,$ by using a diagonal process and passing to a further subsequence if necessary, we can assume that
\begin{align}\label{4.22}
(v^{n_j}(t_{n_j}-J),T^{n_j}(t_{n_j}-J))\rightharpoonup (v_J,T_J)\,\,\textit{weakly in}\,\,V,\,\,\forall\,\,J\in\mathbb{N}
\end{align}
with $(v_J,T_J)\in B_0.$ We deduce from Theorem \ref{4.5} that
\begin{align*}
v=\lim_{w-j\rightarrow+\infty} v^{n_j}(t_{n_j})=\lim_{w-j\rightarrow+\infty}S_1(J) (v^{n_j}(t_{n_j}-J),T^{n_j}(t_{n_j}-J))=S_1(J)(v_J,T_J),
\end{align*}
where $\lim_{w-j\rightarrow+\infty}$ is the limit taken in the weak topology of $H_1.$ Thus, we obtain
\begin{align*}
v=S_1(J)(v_J,T_J),\,\,\forall\,\,J\in\mathbb{N}.
\end{align*}
It follows from \eqref{4.20} that
\begin{align}\label{4.23}
\|v\|_2\leq\liminf_{j\rightarrow+\infty} \|(v^{n_j}(t_{n_j})\|_2.
\end{align}
In what follows, we will prove that
\begin{align}\label{4.24}
\limsup_{j\rightarrow+\infty} \|v^{n_j}(t_{n_j})\|_2\leq\|v\|_2.
\end{align}
For $J\in\mathbb{N}$ and $t_n>J,$ we conclude from \eqref{4.19} that
\begin{align}\label{4.25}
\nonumber&\|v^n(t_n)\|^2_2=e^{-\delta J}\|v^n(t_n-J)\|^2_2-2\int_0^Je^{-\delta(J-s)}a(S_1(s)(v^n(t_n-J),T^n(t_n-J)),S_1(s)(v^n(t_n-J),T^n(t_n-J)))\,ds\\
&+2\int_0^Je^{-\delta(J-s)}\left(\int_{\Omega}(\int_0^z\nabla T^n(t_n+s-J)\,dr)\cdot S_1(s)(v^n(t_n-J),T^n(t_n-J))\,dxdydz\right)\,ds.
\end{align}
From \eqref{4.21}-\eqref{4.22} and Theorem \ref{4.5}, we conclude that
\begin{align}\label{4.26}
e^{-\delta J}\|v^{n_j}(t_{n_j}-J)\|^2_2\leq \rho^2e^{-\delta J}
\end{align}
and
\begin{align*}
S_1(\cdot)(v^{n_j}(t_{n_j}-J),T^{n_j}(t_{n_j}-J))\rightharpoonup S_1(\cdot)(v_J,T_J)\,\,\textit{weakly in}\,\,L^2(0,J;V_1),
\end{align*}
which entails that
\begin{align}\label{4.27}
\nonumber&\int_0^Je^{-\delta(J-s)}a(S_1(s)(v_J,T_J),S_1(s)(v_J,T_J))\,ds\\
\leq&\liminf_{j\rightarrow+\infty}\int_0^Je^{-\delta(J-s)}a(S_1(s)(v^{n_j}(t_{n_j}-J),T^{n_j}(t_{n_j}-J)),S_1(s)(v^{n_j}(t_{n_j}-J),T^{n_j}(t_{n_j}-J)))\,ds.
\end{align}
Thanks to the asymptotical compactness and continuity of the semigroup $\{S_2(t)\}_{t \geq 0}$ in $V_2,$ we obtain
\begin{align}\label{4.28}
\nonumber&\lim_{j\rightarrow+\infty}\int_0^Je^{-\delta(J-s)}\left(\int_{\Omega}(\int_0^z\nabla T^{n_j}(t_{n_j}+s-J)\,dr)\cdot S_1(s)(v^{n_j}(t_{n_j}-J),T^{n_j}(t_{n_j}-J))\,dxdydz)\right)\,ds\\
=&\int_0^Je^{-\delta(J-s)}\left(\int_{\Omega}(\int_0^z\nabla S_2(s)(v_J,T_J)\,dr)\cdot S_1(s)(v_J,T_J)\,dxdydz\right)\,ds.
\end{align}
From \eqref{4.25}-\eqref{4.28}, we infer that
\begin{align}\label{4.29}
\nonumber&\limsup_{j\rightarrow+\infty}\|v^{n_j}(t_{n_j})\|^2_2\leq e^{-\delta J}\rho^2-2\int_0^Je^{-\delta(J-s)}a(S_1(s)(v_J,T_J),S_1(s)(v_J,T_J))\,ds\\
&+2\int_0^Je^{-\delta(J-s)}\left(\int_{\Omega}(\int_0^z\nabla S_2(s)(v_J,T_J)\,dr)\cdot S_1(s)(v_J,T_J)\,dxdydz\right)\,ds.
\end{align}
Thanks to
\begin{align}\label{4.30}
\nonumber&\|v\|^2_2=\|S_1(J)(v_J,T_J)\|^2_2\\
\nonumber=&e^{-\delta J}\|v_J\|^2_2-2\int_0^Je^{-\delta(J-s)}a(S_1(s)(v_J,T_J),S_1(s)(v_J,T_J))\,ds\\
+&2\int_0^Je^{-\delta(J-s)}\left(\int_{\Omega}(\int_0^z\nabla S_2(s)(v_J,T_J)\,dr)\cdot S_1(s)(v_J,T_J)\,dxdydz\right)\,ds,
\end{align}
it follows from \eqref{4.29}-\eqref{4.30} that
\begin{align}\label{4.31}
\nonumber\limsup_{j\rightarrow+\infty}\|v^{n_j}(t_{n_j})\|^2_2\leq &e^{-\delta J}(\rho^2-\|v_J\|^2_2)+\|v\|^2_2\\
\leq &e^{-\delta J}\rho^2+\|v\|^2_2,\,\,\,\forall\,\,J\in\mathbb{N}.
\end{align}
Let $J\rightarrow+\infty,$ we obtain
\begin{align}\label{4.31}
\limsup_{j\rightarrow+\infty}\|v^{n_j}(t_{n_j})\|^2_2\leq \|v\|^2_2.
\end{align}
Combining \eqref{4.23} with \eqref{4.31}, we conclude that
\begin{align}\label{4.31}
\lim_{j\rightarrow+\infty}\|v^{n_j}(t_{n_j})-v\|_2=0.
\end{align}

Due to $\|v\|\leq C\|v\|_2^\frac{1}{2}\|v\|_{(H^2(\Omega))^2\cap V_1}^\frac{1}{2},$ we conclude that the semigroup $\{S_1(t)\}_{t \geq 0}$ is asymptotically compact in $V_1.$ Therefore, the semigroup $\{S(t)\}_{t \geq 0}$
generated by problem \eqref{2.4}-\eqref{2.6} is asymptotically compact in $V.$
The proof of Theorem \ref{4.13} is completed.\\
\hfill\qed

From Theorem \ref{3.1.1}, we know that the semigroup generated by problem \eqref{2.4}-\eqref{2.6} is continuous in $V.$ Furthermore, combining Theorem \ref{4.4} with Theorem \ref{4.13}, we immediately obtain the following result.
\begin{theorem}
Assume that $Q \in H^1(\Omega).$ Then the semigroup $\{S(t)\}_{t \geq 0}$
corresponding to problem \eqref{2.4}-\eqref{2.6} possesses a global
attractor $\mathcal{A}$ in $V.$
\end{theorem}

\section*{Acknowledgement}
 This work was supported by the National Science Foundation of China Grant (11401459) and the Natural Science
Foundation of Shaanxi Province (2015JM1010).

\bibliographystyle{elsarticle-template-num}
\bibliography{BIB}
\end{document}